\documentclass[12pt, journal, draftclsnofoot, onecolumn]{IEEEtran}
\ifCLASSOPTIONcompsoc
\else
\fi

\ifCLASSINFOpdf
\else
\fi
\usepackage[parfill]{parskip}    
\usepackage{amsmath,amssymb,latexsym,enumerate,mathrsfs}
\usepackage{graphicx}
\usepackage{hyperref}
\usepackage{color,soul}
\definecolor{lightblue}{rgb}{.90,.95,1}
\sethlcolor{yellow}

\usepackage{array}
\usepackage{times}
\DeclareMathOperator*{\arginf}{arg\,inf}

\DeclareMathOperator*{\essinf}{ess\,inf}

\newtheorem{theorem}{Theorem}[section]

\newtheorem{lemma}[theorem]{Lemma}
\newtheorem{proposition}[theorem]{Proposition}

\newtheorem{property}{Property}
\newtheorem{definition}[theorem]{Definition}
\newtheorem{assumption}[theorem]{Assumption}
\newtheorem{remark}[theorem]{Remark}

\def\etal{\mbox{et al.}}

\begin{document}
%
\title{On the hierarchical risk-averse control problems for diffusion processes}

\author{Getachew K. Befekadu, \IEEEmembership{} Alexander Veremyev and Eduardo L. Pasiliao \IEEEmembership{}
\thanks{Version -- January, 2018 }
\IEEEcompsocitemizethanks{\IEEEcompsocthanksitem G. K. Befekadu is with the National Research Council, 
Air Force Research Laboratory \& Department of Industrial System Engineering, University of Florida - REEF, 1350 N. Poquito Rd, Shalimar, FL 32579, USA. \protect\\
E-mail: gbefekadu@ufl.edu
\IEEEcompsocthanksitem A. Veremyev is with the Department of Industrial Engineering and Management Systems, University of Central
Florida, 12800 Pegasus Dr., Orlando, FL 32816, USA.\protect\\
E-mail: alexander.veremyev@ucf.edu
\IEEEcompsocthanksitem E. L. Pasiliao is with the Munitions Directorate, Air Force Research Laboratory, 101 West Eglin Blvd, Eglin AFB, FL 32542, USA.\protect\\
E-mail: pasiliao@eglin.af.mil}}

\markboth{}%
{Shell \MakeLowercase{\textit{et al.}}: Bare Advanced Demo of IEEEtran.cls for Journals}
\IEEEcompsoctitleabstractindextext{%
\begin{abstract}
In this paper, we consider a risk-averse control problem for diffusion processes, in which there is a partition of the admissible control strategy into two decision-making groups (namely, the {\it leader} and {\it follower}) with different cost functionals and risk-averse satisfactions. Our approach, based on a hierarchical optimization framework, requires that a certain level of risk-averse satisfaction be achieved for the {\it leader} as a priority over that of the {\it follower's} risk-averseness. In particular, we formulate such a risk-averse control problem involving a family of time-consistent dynamic convex risk measures induced by conditional $g$-expectations (i.e., filtration-consistent nonlinear expectations associated with the generators of certain backward stochastic differential equations). Moreover, under suitable conditions, we establish the existence of optimal risk-averse solutions, in the sense of viscosity solutions, for the corresponding risk-averse dynamic programming equations. Finally, we briefly comment on the implication of our results.

\end{abstract}

}

\maketitle

\IEEEdisplaynotcompsoctitleabstractindextext

%
\IEEEpeerreviewmaketitle

\section{Introduction}	 \label{S1}
Let $\bigl(\Omega, \mathcal{F},\{\mathcal{F}_t \}_{t \ge 0}, \mathbb{P}\bigr)$ be a probability space, and let $\{B_t\}_{t \ge 0}$ be a $d$-dimensional standard Brownian motion, whose natural filtration, augmented by all $\mathbb{P}$-null sets, is denoted by $\{\mathcal{F}_t\}_{t \ge 0}$, so that it satisfies the {\it usual hypotheses} (e.g., see \cite{Pro90} or \cite{GikS72}). \,We consider the following controlled-diffusion process over a given finite-time horizon $T>0$
\begin{align}
d X_t = f\bigl(t, X_t, u_t\bigr) dt + \sigma\bigl(t, X_t, u_t\bigr)dB_t, \quad X_0=x, \quad 0 \le t \le T, \label{Eq1.1} 
\end{align}
where
\begin{itemize}
\item $X_{\cdot}$ is an $\mathbb{R}^{d}$-valued diffusion process,
\item $u_{\cdot}$ is a $U$-valued measurable admissible control process (where $U$ is a compact set in $\mathbb{R}^{d}$) such that for all $t > s$, $(B_t - B_s)$ is independent of $u_r$ for $r \le s$ (nonanticipativity condition) and
\begin{align*}
\mathbb{E} \int_{s}^{t} \vert u_{\tau}\vert^2 d\tau < \infty, \quad \forall t \ge s, 
\end{align*}
\item the function $f \colon [0, T] \times \mathbb{R}^d \times U \rightarrow \mathbb{R}^{d}$ is uniformly Lipschitz, with bounded first derivative, and
\item $\sigma \colon [0, T] \times \mathbb{R}^{d} \times U \rightarrow \mathbb{R}^{d \times d}$ is Lipschitz with the least eigenvalue of $\sigma\,\sigma^T$ uniformly 
bounded away from zero, i.e., 
\begin{align*}
 y^T \sigma(t, x, u)\,\sigma^T(t, x, u) y  \ge \lambda I_{d} \Vert y \Vert^2,  \,\, \forall (x,u) \in \mathbb{R}^{d} \times U, \,\, y \in \mathbb{R}^{d}, \,\,  \forall t \in [0, T],  
\end{align*}
for some $\lambda > 0$.
\end{itemize}

In this paper, we specifically consider a hierarchical risk-averse control problem for the above controlled-diffusion process in \eqref{Eq1.1}, in which there is a partition of the admissible control strategy into two decision making groups (i.e., progressively measurable strategies corresponding to the {\it leader} and {\it follower} -- where such notions are used in the Stackelberg's optimization \cite{VonSta34}) with different cost functionals and risk-averse satisfactions. In particular, we partition the control subdomain $U$ into two open sets $V$ and $W$ that are compatible with the strategy subspaces of the {\it leader} and the {\it follower}, respectively. {More specifically, the risk-averse strategy for the {\it leader} $v_{\cdot}$ is a $V$-valued measurable control process from the set} $\mathcal{V}_{[0,T]}$ {with}
\begin{align}
\mathcal{V}_{[0,T]} \triangleq \Bigl\{v\colon [0,T] \times \Omega \rightarrow \, V \,\Bigl\vert  v \,\, \text{is an} \,\, \bigl\{\mathcal{F}_t\bigr\}_{t\ge 0}\text{- adapted} \,\,  \text{and}\,\, \mathbb{E} \int_{0}^{T} \vert v_t\vert^2 dt < \infty \Bigr\} \label{Eq1.3} 
\end{align}
{and the risk-averse strategy for the} {\it follower} {$w_{\cdot}$ is a $W$-valued measurable control process from the set $\mathcal{W}_{[0,T]}$ with}
\begin{align}
\mathcal{W}_{[0,T]} \triangleq \Bigl\{w\colon [0,T] \times \Omega \rightarrow  \,W \,\Bigl\vert w \,\, \text{is an} \,\, \bigl\{\mathcal{F}_t\bigr\}_{t\ge 0}\text{- adapted} \,\,  \text{and}\,\, \mathbb{E} \int_{0}^{T} \vert w_t\vert^2 dt < \infty \Bigr\}, \label{Eq1.4} 
\end{align}
{respectively.}

Furthermore, we consider the following two cost functionals that provide information about the accumulated risk-costs on the time interval $[0, T]$ w.r.t. the strategies of the {\it leader} and that of the {\it follower}, i.e., 
\begin{align}
 &\text{\it leader's accumulated risk-cost:} \notag \\
 &\quad  \xi_{0,T}^1(v,w) = \int_0^T c_1\bigl(t, X_t, v_t\bigr) dt + \Psi_1(X_T) \label{Eq1.5} 
\end{align}
and
\begin{align}
&\text{\it follower's accumulated risk-cost:} \notag \\
& \quad  \xi_{0,T}^2(v,w) = \int_0^T c_2\bigl(t, X_t, w_t\bigr) dt + \Psi_2(X_T), \label{Eq1.6} 
\end{align}
where the cost-rate functionals $c_1 \colon [0,T] \times \mathbb{R}^d \times V \rightarrow \mathbb{R}$ and $c_2 \colon [0,T] \times \mathbb{R}^d \times W \rightarrow \mathbb{R}$ are measurable functions; and {$\Psi_i\colon \mathbb{R}^d \rightarrow \mathbb{R}$, with $i=1, 2$, are the final-stage risk-costs (that are associated with risk-averse satisfaction levels) and they are assumed to be measurable functions. Here, we also assume that $f$, $\sigma$, $c_i$ and $\Psi_i$, $i=1,2$, satisfy the following growth conditions}
\begin{align}
\bigl\vert f\bigl(t, x, (v, w)) \bigr\vert &+ \bigl\vert \sigma\bigl(t, x, (v,w) \bigr) \bigr\vert + \bigl\vert c_1\bigl(t, x, v\bigr) \bigr\vert + \bigl\vert \Psi_1\bigl(x\bigr) \bigr\vert  \notag \\
 & \quad \le K \bigl(1 + \bigl\vert x \bigr\vert^p + \bigl\vert v \bigr\vert + \bigl\vert w \bigr\vert \bigr) \label{Eq1.7}
\end{align}
and
\begin{align}
\bigl\vert f\bigl(t, x, (v, w)) \bigr\vert &+ \bigl\vert \sigma\bigl(t, x, (v,w) \bigr) \bigr\vert + \bigl\vert c_2\bigl(t, x, w\bigr) \bigr\vert + \bigl\vert \Psi_2\bigl(x\bigr) \bigr\vert \notag \\
 & \quad \le K \bigl(1 + \bigl\vert x \bigr\vert^p + \bigl\vert v \bigr\vert + \bigl\vert w \bigr\vert \bigr), \label{Eq1.8}
 \end{align}
 {for all $\bigl(t, x, (v,w) \bigr) \in [0,T] \times \mathbb{R}^{d} \times (V \times W)$, $p \ge 1$ and for some constant $K > 0$}. Note that the corresponding solution $X_t$ in \eqref{Eq1.1} (i.e., $X_t = X_t^{0,x;u}$ with $u_{\cdot} \equiv (v_{\cdot}, w_{\cdot}) \in \mathcal{V}_{[0,T]} \times\mathcal{W}_{[0,T]}$) depends on the admissible risk-averse strategies of the {\it leader} and that of the {\it follower}; and, moreover, it also depends on the initial condition $X_0=x$.  As a result of this, for any time-interval $[t, T]$, with $t \in [0, T]$, the accumulated risk-costs $ \xi_{t,T}^1$ and $\xi_{t,T}^2$ depend on the risk-averse strategies $(v_{\cdot}, w_{\cdot}) \in \mathcal{V}_{[t,T]} \times\mathcal{W}_{[t,T]}$. 

Next, let us introduce the following measurable spaces that will be useful later in the \,paper.
\begin{itemize}
\item $L^2\bigl(\Omega, \mathcal{F}_t, \mathbb{P}; \mathbb{R}^{d} \bigr)$ is the set of $\mathbb{R}^{d}$-valued $\mathcal{F}_t$-measurable random variables $\xi$ such that $\bigl\Vert \xi \bigr\Vert^2 = \mathbb{E}\bigl\{\bigl\vert \xi \bigr\vert^2  \bigr\}< \infty$;
\item $L^{\infty}\bigl(\Omega, \mathcal{F}_t, \mathbb{P}\bigr)$ is the set of $\mathbb{R}$-valued $\mathcal{F}_t$-measurable random variables $\xi$ such that $\bigl\Vert \xi \bigr\Vert = \essinf \bigl\vert \xi \bigr \vert < \infty$;
\item $\mathcal{S}^2\bigl(t, T; \mathbb{R}^{d} \bigr)$ is the set of $\mathbb{R}^{d}$-valued adapted processes $\bigl (\varphi_s\bigr)_{t \le s \le T}$ on $\Omega \times [t, T]$ such that $\bigl\Vert \varphi \bigr\Vert_{[t, T]}^2 = \mathbb{E}\bigl\{\sup_{t \le s \le T} \bigl\vert \varphi_s \bigr\vert^2  \bigr\}< \infty$;
\item $\mathcal{H}^2\bigl(t, T; \mathbb{R}^{d} \bigr)$ is the set of $\mathbb{R}^{d}$-valued progressively measurable processes $\bigl (\varphi_s\bigr)_{t \le s \le T}$ such that $\bigl\Vert \varphi \bigr\Vert_{[t, T]}^2 = \mathbb{E}\bigl\{ \int_t^T \bigl\vert \varphi_s \bigr \vert^2 ds  \bigr\}< \infty$.
\end{itemize}

On the same probability space $\bigl(\Omega, \mathcal{F},\{\mathcal{F}_t \}_{t \ge 0}, \mathbb{P}\bigr)$, we consider the following one-dimensional backward stochastic differential equation (BSDE) 
\begin{align}
- d Y_t = g\bigl(t, Y_t, Z_t\bigr) dt - Z_tdB_t, \quad Y_T=\xi, \label{Eq1.9}
\end{align}
where the terminal value $Y_T=\xi$ belongs to $L^2\bigl(\Omega, \mathcal{F}_T, \mathbb{P}; \mathbb{R}\bigr)$ and the generator function $g \colon \Omega \times [0, T] \times \mathbb{R} \times \mathbb{R}^{d} \rightarrow \mathbb{R}$, with property that $\bigl(g\bigl(t, y, z\bigr)\bigr)_{0 \le t \le T}$ is progressively measurable for each $(y, z) \in \mathbb{R} \times \mathbb{R}^{d}$. We also assume that $g$ satisfies the following assumption.

\begin{assumption} \label{AS1} 
\begin{enumerate} 
\item[{\rm (A1)}] $g$ is Lipschitz in $(y, z)$, i.e., there exists a constant $K > 0$ such that, $\mathbb{P}$-a.s., for any $t \in [0, T]$, $y_1, y_2 \in \mathbb{R}$ and $z_1, z_2 \in \mathbb{R}^d$ 
\begin{align*}
\bigl\vert g\bigl(t, y_1, z_1\bigr) - g\bigl(t, y_2, z_2\bigr) \bigr\vert \le K \bigl(\bigl\vert y_1 - y_2 \bigr\vert + \bigl\Vert z_1 - z_2 \bigr\Vert\bigr).
\end{align*}
\item[{\rm (A2)}] $g\bigl(t, 0, 0\bigr) \in \mathcal{H}^2\bigl(t, T; \mathbb{R} \bigr)$.
\item[{\rm (A3)}] $\mathbb{P}$-a.s., for all $t \in [0, T]$ and $y \in \mathbb{R}$, $g\bigl(t, y, 0\bigr) = 0$.
\end{enumerate}
\end{assumption}
Then, we state the following lemma, which is used to establish the existence of a unique adapted solution (e.g., see \cite{ParP90} for additional discussions).
\begin{lemma} \label{L1}
Suppose that Assumption~\ref{AS1} holds. Then, for any $\xi \in L^2\bigl(\Omega, \mathcal{F}_T, \mathbb{P}; \mathbb{R}\bigr)$, the BSDE in \eqref{Eq1.9}, with terminal condition $Y_T=\xi$, i.e.,
\begin{align}
 Y_t = \xi + \int_t^T g\bigl(s, Y_s, Z_s\bigr) ds - \int_t^TZ_s dB_s, \quad 0 \le t \le T  \label{Eq1.10}
\end{align}
has a unique adapted solution
\begin{align}
 \bigl(Y_t^{T,g,\xi}, Z_t^{T,g,\xi}\bigr)_{0 \le t \le T} \in  \mathcal{S}^2\bigl(0, T; \mathbb{R} \bigr) \times  \mathcal{H}^2\bigl(0, T; \mathbb{R}^{d} \bigr). \label{Eq1.11}
\end{align}
\end{lemma}
Moreover, we recall the following comparison result that will be useful later (e.g., see \cite{ParT99}).
\begin{theorem} \label{T1}
Given two generators $g_1$ and $g_2$ satisfying Assumption~\ref{AS1} and two terminal conditions $\xi_1, \xi_2 \in L^2\bigl(\Omega, \mathcal{F}_T, \mathbb{P}; \mathbb{R}\bigr)$. Let $\bigl(Y_t^1, Z_t^1\bigr)$ and $\bigl(Y_t^2, Z_t^2\bigr)$ be the solution pairs corresponding to $\bigl(\xi_1, g_1\bigr)$ and $\bigl(\xi_2, g_2\bigr)$, respectively. Then, we have 
\begin{enumerate} [(i)]
\item Monotonicity: If $\xi_1 \ge \xi_2$ and $g_1 \ge g_2$, $\mathbb{P}$-a.s., then $Y_t^1 \ge Y_t^2$, $\mathbb{P}$-a.s., for all $t \in [0, T]$;
\item Strictly Monotonicity: In addition to ($i$) above, if we assume that $\mathbb{P}\bigl(\xi_1 > \xi_2\bigr) > 0$, then $\mathbb{P}\bigl(Y_t^1 > Y_t^2\bigr) > 0$, for all $t \in [0, T]$.
\end{enumerate}
\end{theorem}
In the following, we give a definition for a dynamic risk measure that is associated with the generator of the BSDE in \eqref{Eq1.9}.
\begin{definition} \label{Df1} 
For any $\xi \in L^2\bigl(\Omega, \mathcal{F}_T, \mathbb{P}; \mathbb{R}\bigr)$, let $\bigl(Y_t^{T,g,\xi}, Z_t^{T,g,\xi}\bigr)_{0 \le t \le T} \in  \mathcal{S}^2\bigl(0, T; \mathbb{R} \bigr) \times  \mathcal{H}^2\bigl(0, T; \mathbb{R}^{d} \bigr)$ be the unique solution to the BSDE in \eqref{Eq1.9} with terminal condition $Y_T=\xi$. Then, we define the dynamic risk measure $\rho_{t,T}^g$ of $\xi$ by\footnote{{Note that, for any $t \in [0,T]$, the conditional $g$-expectation of $\xi \in L^2\bigl(\Omega, \mathcal{F}_T, \mathbb{P}; \mathbb{R}\bigr)$ w.r.t. $\mathcal{F}_t$ (denoted by $\mathcal{E}_g\bigl[\xi \vert \mathcal{F}_t\bigr]$ (e.g., see} \cite{Ros06} or \cite{Pen04}; cf. \cite{Bio09})) {is also defined by}
\begin{align*}
\mathcal{E}_g\bigl[\xi \vert \mathcal{F}_t\bigr] \triangleq Y_t^{T,g,\xi}.
\end{align*}
{Moreover, if $\xi$ is $\mathcal{F}_t$-measurable, then we have $\mathcal{E}_g\bigl[\xi \vert \mathcal{F}_t\bigr] = \xi$.} }
\begin{align}
\rho_{t,T}^g \bigl[\xi \bigr] \triangleq Y_t^{T,g,\xi}.  \label{Eq1.12}
\end{align}
\end{definition}

Note that such a risk measure is widely used for evaluating the risk of stochastic processes or uncertain outcomes, and assists with stipulating minimum interventions required by financial institutions for risk management (e.g., see \cite{ArtDEH99}, \cite{Pen04}, \cite{ElKPQ97}, \cite{FolS02}, \cite{DetS05} or \cite{CorHMP02} for related discussions). In Section~\ref{S2}, we use a family of dynamic risk measures associated {that takes into account} the {\it leader's} and {\it follower's} cost functionals and risk-averse satisfactions; and we further provide a hierarchical framework for the risk-averse control problem.

Moreover, if the generator function $g$ satisfies Assumption~\ref{AS1}, then a family of time-consistent {dynamic risk measures} $\bigl\{\rho_{t,T}^g\bigr\}_{t \in [0,T]}$ has the following properties (e.g., see \cite{Pen04} for additional discussions). 
\begin{property} \label{Pr1} 
\begin{enumerate} 
\item[{\rm (P1)}] {\it Normalization}: $\rho_{t,T}^g\bigl[0\bigr] = 0$ for $t \in [0, T]$,
\item[{\rm (P2)}] {\it Monotonicity}:  For $\xi_1, \xi_2 \in L^2\bigl(\Omega, \mathcal{F}_T, \mathbb{P}; \mathbb{R}\bigr)$ such that $\xi_1 \ge \xi_2$ $\mathbb{P}$-a.s., then 
\begin{align*}
\rho_{t,T}^g\bigl[\xi_1\bigr] \ge \rho_{t,T}^g\bigl[\xi_2\bigr], \quad \mathbb{P}{\text-a.s.},
\end{align*}
\item[{\rm (P3)}] {\it Translation Invariance}: For all $\xi \in L^2\bigl(\Omega, \mathcal{F}_T, \mathbb{P}; \mathbb{R}\bigr)$ and $\nu \in L^2\bigl(\Omega, \mathcal{F}_t, \mathbb{P}; \mathbb{R}\bigr)$
\begin{align*}
\rho_{t,T}^g\bigl[\xi + \nu\bigr] = \rho_{t,T}^g\bigl[\xi\bigr] + \nu,
\end{align*}
\item[{\rm (P4)}] {\it Convexity}: If $g$ is a convex function for every fixed $(t, \omega) \in [0, T] \times \Omega$, then for all $\xi_1, \xi_2 \in L^2\bigl(\Omega, \mathcal{F}_T, \mathbb{P}; \mathbb{R}\bigr)$ and for all $\lambda \in L^{\infty}\bigl(\Omega, \mathcal{F}_t, \mathbb{P}; \mathbb{R}\bigr)$ such that $0 \le \lambda \le 1$
\begin{align*}
\rho_{t,T}^g\bigl[\lambda \xi_1 + (1-\lambda)\xi_2 \bigr] \le \lambda \rho_{t,T}^g\bigl[\xi_1\bigr] + (1- \lambda) \rho_{t,T}^g\bigl[\xi_1\bigr],
\end{align*}
\item[{\rm (P5)}] {\it Positive Homogeneity}: {If $g$ is positive-homogeneous,} then, for all $\xi \in L^2\bigl(\Omega, \mathcal{F}_T, \mathbb{P}; \mathbb{R}\bigr)$ and for all $\lambda \in L^{\infty}\bigl(\Omega, \mathcal{F}_t, \mathbb{P}; \mathbb{R}\bigr)$ such that $\lambda > 0$
\begin{align*}
\rho_{t,T}^g\bigl[\lambda \xi \bigr] = \lambda \rho_{t,T}^g\bigl[\xi\bigr].
\end{align*}
\end{enumerate}
\end{property}

\begin{remark}
Note that, since the seminal work of Artzner \etal~\cite{ArtDEH99}, there have been studies on axiomatic dynamic risk measures, coherency and consistency in the literature (e.g., see \cite{DetS05}, \cite{Pen04}, \cite{Ros06}, \cite{FolS02} or \cite{CorHMP02}). Particularly relevant for us is {a family of time-consistent dynamic convex risk measures, based on the conditional $g$-expectations associated with} the generators of BSDEs in \eqref{Eq2.8} and \eqref{Eq2.9}, that satisfy the above properties {\rm (P1)}--{\rm (P5)}.
\end{remark}

Here, it is worth mentioning that some interesting studies on the dynamic risk measures, {based on the conditional $g$-expectations}, have been reported in the literature (e,g. see \cite{Pen04}, \cite{CorHMP02} and \cite{Ros06} for establishing connection between the risk measures and the generator of BSDE; and see also \cite{RSta10} for characterizing the generator of BSDE according to different risk measures). Recently, the authors in \cite{Rus10} and \cite{BefVP16} (see also \cite{BefVBP17}) have provided interesting results on the risk-averse control problem for Markov decision processes, in discrete-time setting, and, respectively, for controlled-diffusion processes in continuous-time framework. {Note that the rationale behind our framework, which follows in some sense the settings of these papers, is to show how a hierarchical optimization framework can be systematically used to obtain optimal risk-averse strategies for such controlled-diffusion processes.}\footnote{In this paper, our intent is to provide a theoretical framework, rather than considering a specific numerical problem or application.}

The remainder of this paper is organized as follows. In Section~\ref{S2}, using the basic remarks made in Section~\ref{S1}, we state the risk-averse control problem for the controlled-diffusion process. In Section~\ref{S3}, we present our main results -- where we introduce a framework under which the {\it follower} is required to {\it respond optimally} to the risk-averse strategy of the {\it leader} so as to achieve {an overall optimal risk-averseness. In this section}, we also establish the existence of optimal risk-averse solutions, in the sense of viscosity solutions, to the associated risk-averse dynamic programming equations. Finally, Section~\ref{S4} provides further remarks.

\section{The hierarchical risk-averse control problem formulation} \label{S2}
{In order to make our formulation more precise, for any $(t, x) \in [0, T] \times \mathbb{R}^d$, we consider the following forward SDE with an initial condition $X_t^{t,x;u} = x$}
\begin{align}
d X_s^{t,x;u} = & f\bigl(t, X_s^{t,x; u}, (v_s, w_s)\bigr) ds + \sigma\bigl(s, X_s^{t,x; u}, (v_s, w_s)\bigr)dB_s, \,\,  t \le s \le T,  \label{Eq2.1}
\end{align}
where $v_{\cdot}$ and $w_{\cdot}$ are $V$- and $W$-valued measurable control processes, respectively. Moreover, we introduce the following two risk-value functions w.r.t. the strategies of the {\it leader} and that of the {\it follower}, i.e.,
\begin{align}
 &\text{\it leader:}  \quad V_1^{v}\bigl(t, x\bigr) = \rho_{t, T}^{g_1} \bigl[\xi_{t,T}^1\bigl(v, \hat{w}\bigr)\bigr],  \,\,\notag \\
 &\quad \quad  \text{such that} \notag \\
 & \quad  \quad  \hat{w} \in \Bigl\{\tilde{w}_{\cdot} \in \mathcal{W}_{[t,T]} \,\Bigl\vert\, \rho_{t, T}^{g_2} \bigl[\xi_{t,T}^2\bigl(v, \tilde{w}\bigr)\bigl] \le \rho_{t, T}^{g_2} \bigl[\xi_{t,T}^2\bigl(v, w\bigr)\bigl], \notag \\
   &  \hspace{0.9in}   \forall w_{\cdot} \in \mathcal{W}_{[t,T]}, \,\, v_{\cdot} \in \mathcal{V}_{[t,T]}, \,\, \mathbb{P}-a.s. \Bigr\}, \label{Eq2.2}
\end{align}
where
\begin{align}
\xi_{t,T}^1\bigl(v, \hat{w}\bigr) = \int_t^T c_1\bigl(s, X_s^{t,x;u}, v_s\bigr) ds + \Psi_1(X_T^{t,x;u}), \label{Eq2.3}
\end{align}
{with $u=(v, \hat{w})$;} and similarly
\begin{align}
\text{\it follower:} \quad V_2^{w}\bigl(t, x\bigr) = \rho_{t, T}^{g_2} \bigl[\xi_{t,T}^2\bigl(\hat{v},w\bigr)\bigr], \label{Eq2.4}
\end{align}
where {$\hat{v}_{\cdot} \in \mathcal{V}_{[t,T]}$ is assumed to be known to the} {\it follower} {and }
\begin{align}
\xi_{t,T}^2 \bigl(\hat{v},w\bigr) = \int_t^T c_2\bigl(s, X_s^{t,x;u}, w_s\bigr) ds + \Psi_2(X_T^{t,x;u}), \label{Eq2.5}
\end{align}
{with $u=(\hat{v}, w)$}.

Note that we can express the above risk-value functions in \eqref{Eq2.2} and \eqref{Eq2.4} as follow
\begin{align}
 V_1^{v}\bigl(t, x\bigr) &= \xi_{t,T}^1\bigl(v,\hat{w}\bigr) + \int_t^T g_1\bigl(s, Y_s^{1;t,x;u}, Z_s^{1;t,x;u}\bigr) ds - \int_t^T Z_s^{1;t,x;u} dB_s \notag \\
                                                                      &= \Psi_1(X_T^{t,x;u}) + \int_t^T \bigl\{c_1\bigl(s, X_s^{t,x;u}, v_s\bigr)  + g_1\bigl(s, Y_s^{1;t,x;u}, Z_s^{1;t,x;u} \bigr) \bigr\}ds - \int_t^T Z_s^{1;t,x;u}dB_s, \label{Eq2.6}
\end{align}
{with $u=(v,\hat{w})$ and} $\hat{w} \in \bigl\{\tilde{w}_{\cdot} \in \mathcal{W}_{[t,T]} \,\Bigl\vert\, \rho_{t, T}^{g_2} \bigl[\xi_{t,T}^2\bigl(v, \tilde{w}\bigr)\bigl] \le \rho_{t, T}^{g_2} \bigl[\xi_{t,T}^2\bigl(v, w\bigr)\bigl], \,\, \forall w_{\cdot} \in \mathcal{W}_{[t,T]}, \,\, v_{\cdot} \in \mathcal{V}_{[t,T]} \bigr\}$, and 
\begin{align}
 V_2^{w}\bigl(t, x\bigr) &= \xi_{t,T}^2\bigl(\hat{v},w\bigr) + \int_t^T g_2\bigl(s, Y_s^{2;t,x;u}, Z_s^{2;t,x;u}\bigr) ds  - \int_t^T Z_s^{2;t,x;u} dB_s \notag \\
                                                                      &= \Psi_2(X_T^{t,x;u}) + \int_t^T \bigl\{c_2\bigl(s, X_s^{t,x;u}, w_s\bigr)  + g_2\bigl(s, Y_s^{2;t,x;u}, Z_s^{2;t,x;u}\bigr) \bigr\}ds - \int_t^T Z_s^{2;t,x;u}dB_s, \label{Eq2.7} 
\end{align}
with $u=(\hat{v}, w)$ -- where the generators $g_1$ and $g_2$ are assumed to satisfy Assumption~\ref{AS1}. Further, noting the conditions in \eqref{Eq1.7} and \eqref{Eq1.8}, then $\bigl(Y_s^{1;t,x;u}, Z_s^{1;t,x;u}\bigr)_{t \le s \le T}$ and $\bigl(Y_s^{2;t,x;u}, Z_s^{2;t,x;u}\bigr)_{t \le s \le T}$ are adapted solutions on $[t, T] \times \Omega$ and belong to $\mathcal{S}^2\bigl(t, T; \mathbb{R} \bigr) \times  \mathcal{H}^2\bigl(t, T; \mathbb{R}^{d} \bigr)$. Equivalently, we can also rewrite \eqref{Eq2.6} and \eqref{Eq2.7} as a family of coupled BSDEs on the probability space $\bigl(\Omega, \mathcal{F}, \mathbb{P}, \{\mathcal{F}_t \}_{t \ge 0})$, i.e.,
\begin{align}
 -d Y_s^{1;t,x;u} &= \bigl\{c_1\bigl(s, X_s^{t,x;u}, v_s\bigr) + g_1\bigl(s, Y_s^{1;t,x;u}, Z_s^{1;t,x;u}\bigr) \bigr\}ds \notag \\
 & \,\, - Z_s^{1;t,x;u} dB_s, \quad s \in [t, T], \quad Y_T^1 = \Psi_1(X_T^{t,x;u}) \label{Eq2.8}
\end{align}
and 
\begin{align}
 -d Y_s^{2;t,x;u} &= \bigl\{c_2\bigl(s, X_s^{t,x;u}, w_s\bigr) + g_2\bigl(s, Y_s^{2;t,x;u}, Z_s^{2;t,x;u}\bigr) \bigr\}ds \notag \\
 & \,\,- Z_s^{2;t,x;u} dB_s,  \quad s \in [t, T], \quad Y_T^2 =\Psi_2(X_T^{t,x;u}). \label{Eq2.9}
\end{align}
 
In what follows, we introduce a hierarchical optimization framework that requires a certain level of risk-averse satisfaction be achieved for the {\it leader} as a priority over that of the {\it follower's} risk-averseness. For example, suppose that the risk-averse strategy of the {\it leader} $\hat{v}_{\cdot} \in \mathcal{V}_{[t,T]}$ is given.\footnote{{Note that such a risk-averse strategy is also assumed to guarantee the} {\it leader} {a certain value of accumulated risk-cost that can not be exceeded, no matter what the} {\it follower's rational-response} (see also Section~\ref{S4}).} Then, the problem of finding an optimal risk-averse strategy for the {\it follower}, i.e., $\hat{w}_{\cdot} \in \mathcal{W}_{[t,T]}$, which minimizes the accumulated risk-cost under $w$ is then reduced to finding an optimal risk-averse solution for
\begin{align}
\inf_{w_{\cdot} \in \mathcal{W}_{[t,T]}}  J_2\bigr[\bigl(\hat{v}, w\bigr)\bigl], \label{Eq2.10}
\end{align}
where 
\begin{align}
J_2\bigr[\bigl(\hat{v}, w\bigr)\bigl] = \rho_{t, T}^{g_2} \bigl[\xi_{t,T}^2\bigl(\hat{v},w\bigr)\bigr]. \label{Eq2.11}
\end{align}
Note that, for a given $\hat{v}_{\cdot} \in \mathcal{V}_{[t,T]} $, if the forward-backward stochastic differential equations (FBSDEs) in \eqref{Eq2.1}, \eqref{Eq2.8} and \eqref{Eq2.9} admit weak solutions, then we have
\begin{align}
 \hat{w} \in S(\hat{v}) \subset \Bigl\{\tilde{w}_{\cdot} \in \mathcal{W}_{[t,T]} \Bigl\vert \, \rho_{t, T}^{g_2} \bigl[\xi_{t,T}^2\bigl(\hat{v}, \tilde{w}\bigr)\bigl] \le \rho_{t, T}^{g_2} \bigl[\xi_{t,T}^2\bigl(\hat{v}, w\bigr)\bigl], 
  \,\, \forall w_{\cdot}  \in \mathcal{W}_{[t,T]},  \,\, \mathbb{P}-a.s. \Bigr\}\label{Eq2.12}
\end{align}
for some (nonanticipating) measurable mapping $S \colon \mathcal{V}_{[t,T]} \rightrightarrows \mathcal{W}_{[t,T]}$. Furthermore, if we substitute $u = (\hat{v},S(\hat{v}))$ into \eqref{Eq2.1}, then the corresponding solution $X_s^{t,x;u}$ depends uniformly on $\hat{v}_{\cdot} \in \mathcal{V}_{[t,T]}$ for $s \in [t, T]$.\footnote{In this paper, for the sake of simplicity, we use the same notation to represent both the set in \eqref{Eq2.12} and an element of that set which is uniquely selectable.} Moreover, the risk-averse control problem (which minimizes the accumulated risk-cost under $v$ w.r..t the {\it leader}) is then reduced to finding an optimal risk-averse solution for
\begin{align}
\inf_{v_{\cdot} \in \mathcal{V}_{[t,T]}}  J_1\bigl[\bigl(v, S(v)\bigr) \bigr], \label{Eq2.13}
\end{align}
where
\begin{align}
J_1\bigr[\bigl(v, S(v)\bigr)\bigl] = \rho_{t, T}^{g_1} \bigl[\xi_{t,T}^1\bigl(v,S(v)\bigr)\bigr]. \label{Eq2.14}
\end{align}

Next, we introduce the definition of admissible hierarchical risk-averse control system $\Sigma_{[t, T]}$ which provides a logical construct for our main results (e.g., see \cite{LIW14}).
\begin{definition} \label{Df2}
For a given finite-time horizon $T>0$, we call $\Sigma_{[t, T]}$ an admissible hierarchical risk-averse control system, if it satisfies the following conditions:
\begin{itemize}
\item $\bigl(\Omega, \mathcal{F},\{\mathcal{F}_t \}_{t \ge 0}, \mathbb{P}\bigr)$ is a complete probability space;
\item $\bigl\{B_s\bigr\}_{s \ge t}$ is a $d$-dimensional standard Brownian motion defined on $\bigl(\Omega, \mathcal{F}, \mathbb{P}\bigr)$ over $[t, T]$ and $\mathcal{F}^t \triangleq \bigl\{\mathcal{F}_s^t\bigr\}_{s \in [t, T]}$, where $\mathcal{F}_s^t = \sigma\bigl\{\bigl(B_s; \,t \le s \le T \bigr)\bigr\}$ is augmented by all $\mathbb{P}$-null sets in $\mathcal{F}$;
\item $v_{\cdot} \colon \Omega \times [s, T]  \rightarrow V$ and $w_{\cdot} \colon \Omega \times [s, T]  \rightarrow W$ are $\bigl\{\mathcal{F}_s^t\bigr\}_{s \ge t}$-adapted processes on $\bigl(\Omega, \mathcal{F}, \mathbb{P}\bigr)$ with 
\begin{align*}
\mathbb{E} \int_{s}^{T} \vert v_{\tau}\vert^2 d \tau < \infty \,\, \text{and}  \,\, \mathbb{E} \int_{s}^{T} \vert w_{\tau}\vert^2 d \tau < \infty, \,\, s \in [t, T];
\end{align*}
\item There exists at least one measurable mapping $S \colon \mathcal{V}_{[t,T]} \rightrightarrows \mathcal{W}_{[t,T]}$ with $w_{\cdot} \in S\bigl(v_{\cdot} \bigr)$ whenever $v_{\cdot} \in \mathcal{V}_{[t,T]}$;
\item For any $x \in \mathbb{R}^d$, the FBSDEs in \eqref{Eq2.1}, \eqref{Eq2.8} and \eqref{Eq2.9} admit a unique solution set \\ $\bigl\{X_{\cdot}^{s,x;u}, (Y_{\cdot}^{1;s,x;u}, Z_{\cdot}^{1;s,x;u}), (Y_{\cdot}^{2;s,x;u}, Z_{\cdot}^{2;s,x;u})\bigr\}$ on $\bigl(\Omega, \mathcal{F}, \mathcal{F}^t, \mathbb{P}\bigr)$ with $u=\bigl(v, S(v)\bigr)$.
\end{itemize}
\end{definition}

Then, for the admissible hierarchical risk-averse control system $\Sigma_{[0, T]}$, we can state the problem of risk-averse control as follow.

{\bf Problem~(P)}. Find a pair of risk-averse strategies $(v_{\cdot}^{\ast}, w_{\cdot}^{\ast}) \in \mathcal{V}_{[0,T]} \times\mathcal{W}_{[0,T]}$ w.r.t. the {\it leader} and that of the {\it follower} such that
\begin{align}
 v_{\cdot}^{\ast} \in \Bigl\{ \arginf J_1\bigr[\bigl(v, w\bigr)\bigl] \Bigl \vert w_{\cdot} \in S(v_{\cdot}) \,\, \& \,\, (v_{\cdot}, S(v_{\cdot})) \,\, \text{restricted to} \,\,\Sigma_{[0, T]} \Bigr\} \subset \mathcal{V}_{[0,T]}  \label{Eq2.15}
\end{align}
and
\begin{align}
  w_{\cdot}^{\ast} \in \Bigl\{ \arginf J_2\bigr[\bigl(v, w\bigr)\bigl] \Bigl \vert  w_{\cdot} \in S(v_{\cdot}) \,\, \& \,\,  (v_{\cdot},S(v_{\cdot})) \,\,  \text{restricted to} \,\Sigma_{[0, T]} \Bigr\} \subset \mathcal{W}_{[0,T]}\label{Eq2.16}
\end{align}
where $S$ is a unique measurable mapping that maps $\mathcal{V}_{[0,T]}$ onto $\mathcal{W}_{[0,T]}$ and, furthermore, the accumulated risk-costs $J_1$ and $J_2$ over the time-interval $[0, T]$ are given  
\begin{align}
J_1\bigr[\bigl(v^{\ast}, w^{\ast}\bigr)\bigl] = \int_0^T c_1\bigl(s, X_s^{0,x;u^{\ast}}, v_s\bigr) ds + \Psi_1(X_T^{0,x;u^{\ast}})  \label{Eq2.17}
\end{align}
and
\begin{align}
J_2\bigr[\bigl(v^{\ast}, w^{\ast}\bigr)\bigl] = \int_0^T c_2\bigl(s, X_s^{0,x;u^{\ast}}, w_s\bigr) ds + \Psi_2(X_T^{0,x;u^{\ast}}),  \label{Eq2.18}
\end{align}
where $X_0^{0,x;u^{\ast}} = x$ and $u^{\ast}=(v^{\ast},w^{\ast})$, with $w_{\cdot}^{\ast} \in S(v_{\cdot}^{\ast})$.

In the following section, we establish the existence of optimal risk-averse solutions, in the sense of viscosity solutions, for the optimization problems in \eqref{Eq2.15} and \eqref{Eq2.16} with restriction to $\Sigma_{[0, T]}$. Note that, for a given $v_{\cdot} \in \mathcal{V}_{[0,T]}$, the optimization problem in \eqref{Eq2.16} has a unique solution on $\mathcal{W}_{[0,T]}$ (see Proposition~\ref{P5}). Furthermore, as we will see later (particularly in Proposition~\ref{P7}), the problem in \eqref{Eq2.15} makes sense if the {\it follower} is involved not only in minimizing his own accumulated risk-cost (in response to the risk-averse strategy of the {\it leader}) but also in minimizing that of the {\it leader}.

\section{Main results} \label{S3}
In this section, we present our main results, {where we establish the existence of optimal risk-averse solutions, in the sense of viscosity solutions, to the associated risk-averse dynamic programming equations (cf. equations~}\eqref{Eq3.14} and \eqref{Eq3.16}).

We now state the following propositions that will be useful for proving our main results later in Subsections~\ref{S3.1} and \ref{S3.2}.
\begin{proposition} \label{P1}
Suppose Assumption~\ref{AS1} together with \eqref{Eq1.7} and \eqref{Eq1.8} hold. Then, for any $(t,x) \in [0, T] \times \mathbb{R}^d$ and for every $(v_{\cdot}, w_{\cdot}) \in \mathcal{V}_{[t,T]} \times\mathcal{W}_{[t,T]}$, the FBSDEs in \eqref{Eq2.1}, \eqref{Eq2.8} and \eqref{Eq2.9} admit unique adapted solutions
\begin{eqnarray}
\left.\begin{array}{c}
X_{\cdot}^{t,x;u} \in  \mathcal{S}^2\bigl(t, T; \mathbb{R} \bigr) ~~~~~~\\
\bigl(Y_{\cdot}^{1;t,x;u}, Z_{\cdot}^{1;t,x;u}\bigr) \in \mathcal{S}^2\bigl(t, T; \mathbb{R} \bigr) \times  \mathcal{H}^2\bigl(t, T; \mathbb{R}^{d} \bigr) \\
\bigl(Y_{\cdot}^{2;t,x;u}, Z_{\cdot}^{2;t,x;u}\bigr) \in \mathcal{S}^2\bigl(t, T; \mathbb{R} \bigr) \times  \mathcal{H}^2\bigl(t, T; \mathbb{R}^{d} \bigr)
\end{array}\right\}  \label{Eq3.1}
\end{eqnarray}
Furthermore, the risk-values w.r.t. the {\it leader} and {\it follower}, i.e., $V_1^{v}\bigl(t, x\bigr)$ and $V_2^{w}\bigl(t, x\bigr)$, are deterministic.
\end{proposition}

\begin{IEEEproof}
Notice that $f$ and $\sigma$ are bounded and Lipschitz continuous w.r.t. $(t,x) \in [0, T] \times \mathbb{R}^d$ and uniformly for $(v, w) \in V \times W$. Then, for any $(t,x) \in [0, T] \times \mathbb{R}^d$ and $u_{\cdot}=(v_{\cdot}, w_{\cdot})$ are progressively measurable processes, there always exists a unique path-wise solution $X_{\cdot}^{t,x;u} \in  \mathcal{S}^2\bigl(t, T; \mathbb{R} \bigr)$ for the forward SDE in \eqref{Eq2.1}. On the other hand, consider the following BSDEs
\begin{align}
 -d \hat{Y}_s^{1;t,x;u} = g_1\bigl(s, Z_s^{1;t,x;u}\bigr) ds - Z_s^{1;t,x;u} dB_s, \label{EqP1.1}
\end{align}
where  
\begin{align*}
\hat{Y}_T^{1;t,x;u} = \int_t^T c_1\bigl(\tau, X_{\tau}^{t,x;u}, v_{\tau}\bigr) d\tau + \Psi_1(X_T^{t,x;u}) 
\end{align*}
and 
\begin{align}
 -d \hat{Y}_s^{2;t,x;u} =  g_2\bigl(s, Z_s^{2;t,x;u}\bigr) ds - Z_s^{2;t,x;u} dB_s, \label{EqP1.2}
\end{align}
where 
\begin{align*}
\hat{Y}_T^{2;t,x;u} = \int_t^T c_2\bigl(\tau, X_{\tau}^{t,x;u}, w_{\tau}\bigr) d\tau + \Psi_2(X_T^{t,x;u}). 
\end{align*}
From Lemma~\ref{L1}, the equations in \eqref{EqP1.1} and \eqref{EqP1.2} admit unique solutions $\bigl(\hat{Y}_{\cdot}^{1;t,x;u}, Z_{\cdot}^{1;t,x;u}\bigr)$ and $\bigl(\hat{Y}_{\cdot}^{2;t,x;u}, Z_{\cdot}^{2;t,x;u}\bigr)$ in $\mathcal{S}^2\bigl(t, T; \mathbb{R} \bigr) \times  \mathcal{H}^2\bigl(t, T; \mathbb{R}^{d} \bigr)$. Furthermore, if we introduce the following
\begin{align*}
Y_s^{1;t,x;u} = \hat{Y}_s^{1;t,x;u} - \int_t^s c_1\bigl(\tau, X_{\tau}^{t,x;u}, v_{\tau}\bigr) d\tau,   \quad s \in [t, T]
\end{align*}
and
\begin{align*}
Y_s^{2;t,x;u} = \hat{Y}_s^{2;t,x;u} - \int_t^s c_2\bigl(\tau, X_{\tau}^{t,x;u}, w_{\tau}\bigr) d\tau,  \quad s \in [t, T].
\end{align*}
Then, the forward SDEs in \eqref{Eq2.8} and \eqref{Eq2.9} hold, with $\bigl(Y_{\cdot}^{1;t,x;u}, Z_{\cdot}^{1;t,x;u}\bigr)$ and $\bigl(Y_{\cdot}^{2;t,x;u}, Z_{\cdot}^{2;t,x;u}\bigr)$, respectively. Moreover, we also observe that $Y_t^{1;t,x;u}$ and $Y_t^{2;t,x;u}$ are deterministic. This completes the proof of Proposition~\ref{P1}. 
\end{IEEEproof}

\begin{proposition} \label{P2}
Let  $(t,x) \in [0, T] \times \mathbb{R}^d$ and $(v_{\cdot}, w_{\cdot}) \in \mathcal{V}_{[t,T]} \times\mathcal{W}_{[t,T]}$ be restricted to $\Sigma_{[t, T]}$ (cf. Definition~\ref{Df2}). Then, for any $r \in [t, T]$ and $\mathbb{R}^d$-valued $\mathcal{F}_r^t$-measurable random variable $\eta$, we have
\begin{align}
 V_1^{v}\bigl(r, \eta\bigr) = \rho_{r, T}^{g_1} \Bigl[\int_r^T c_1\bigl(s, X_s^{r,\eta;u}, v_s\bigr) ds &+ \Psi_1(X_T^{r,\eta;u}) \Bigr], \quad  \mathbb{P}{\text-a.s.}  \label{Eq3.2}
\end{align}
and
\begin{align}
 V_2^{w}\bigl(r, \eta\bigr) = \rho_{r, T}^{g_2} \Bigl[\int_r^T c_2\bigl(s, X_s^{r,\eta;u}, w_s\bigr) ds + \Psi_2(X_T^{r,\eta;u}) \Bigr], \quad \mathbb{P}{\text-a.s.} \label{Eq3.3}
\end{align}
\end{proposition}

\begin{IEEEproof}
For any $r \in [t, T]$, with $t \in [0, T]$, we consider the following probability space $\bigl(\Omega, \mathcal{F}, \mathbb{P}\bigl(\cdot\vert \mathcal{F}_r^t\bigr), \{\mathcal{F}^t\}\bigr)$ and notice that $\eta$ is deterministic under this probability space. Then, for any $s \ge r$, there exist progressively measurable processes $\psi_1$ and $\psi_2$ such that
\begin{align}
 \bigl(v_s(\Omega), w_s(\Omega)\bigl) &= \bigl(\psi_1(\Omega, B_{\cdot \wedge s}(\Omega)), \psi_2(\Omega, B_{\cdot \wedge s}(\Omega)) \bigr), \notag \\  
                                                               &= \bigl(\psi_1(s, \bar{B}_{\cdot \wedge s}(\Omega) + B_{r}(\Omega)), \psi_2(s, \bar{B}_{\cdot \wedge s}(\Omega) + B_{r}(\Omega)) \bigr), \label{EqP2.1}
\end{align}
where $\bar{B}_s = B_s - B_r$ is a standard $d$-dimensional brownian motion. Note that the pairs $\bigl(v_{\cdot}, w_{\cdot}\bigl)$ are $\mathcal{F}_r^t$-adapted processes, then we have the following restriction w.r.t. $\Sigma_{[t, T]}$
\begin{align}
 \bigl(\Omega, \mathcal{F}, \{\mathcal{F}^t\}, \mathbb{P}\bigl(\cdot \vert \mathcal{F}_r^t\bigr)(\omega'), B_{\cdot}, \bigl(v_{\cdot}, w_{\cdot} \bigl)\bigr) \in \Sigma_{[t, T]}, \label{EqP2.2}
\end{align}
where $\omega' \in \Omega'$ such that $\Omega' \in \mathcal{F}$, with $\mathbb{P}(\Omega')=1$. Furthermore, noting Lemma~\ref{L1}, if we work under the probability space $\bigl(\Omega', \mathcal{F}, \mathbb{P}\bigl(\cdot\vert \mathcal{F}_r^t\bigr)\bigr)$, then both statements in \eqref{Eq3.2} and \eqref{Eq3.3} hold $\mathbb{P}$-{\it almost surely}. This completes the proof of Proposition~\ref{P2}. 
\end{IEEEproof}

In what follows, we restrict our discussion when the generators $g_1$ and $g_2$ depend only on $(t, z) \in [0, T] \times \mathbb{R}^d$. Moreover, for $(v, w) \in V \times W$ and any $ \phi(x) \in C_0^{\infty}(\mathbb{R}^d)$, {we introduce the following second-order linear operators}\footnote{{The operators $\mathcal{L}_{t}^{v}$ and $\hat{\mathcal{L}}_{t}^{w}$ are parameterized w.r.t. $v$ and $w$, respectively.}}, {associated with} \eqref{Eq2.1}, as follows
\begin{align}
 \mathcal{L}_{t}^{v} \phi(x)  = \dfrac{1}{2} \operatorname{tr} \Bigl\{a^v(t, x) D_{x}^2 \phi(x)\Bigr\} &+ f(t, x, (v,S(v))) D_{x} \phi(x),  \quad t \in [0, T],  \label{Eq3.4a}
\end{align}
{where $a^v(t, x) = \sigma(t, x,(v, S(v))) \sigma^T(t, x,(v,S(v)))$ and $S$ is the mapping that satisfies} \eqref{Eq2.12}; and
\begin{align}
 \hat{\mathcal{L}}_{t}^{w} \phi(x) = \dfrac{1}{2} \operatorname{tr} \Bigl\{a^w(t, x) D_{x}^2 \phi(x)\Bigr\} &+ f(t, x, (v, w)) D_{x} \phi(x), \quad t \in [0, T],  \label{Eq3.4b}
\end{align}
where $a^w(t, x) = \sigma(t, x,(v,w)) \sigma^T(t, x,(v,w))$, $D_{x}$ and $D_{x}^2$, (with $D_{x}^2 = \bigl({\partial^2 }/{\partial x_i \partial x_j} \bigr)$) are the gradient and the Hessian (w.r.t. the variable $x$), respectively. Furthermore, on the space $C_b^{1,2}([t, T] \times \mathbb{R}^d)$, for any $(t,x) \in [0, T] \times \mathbb{R}^d$, we consider the following coupled Hamilton-Jacobi-Bellman (HJB) partial differential equations
\begin{eqnarray}
\left.\begin{array}{r}
 \dfrac{\partial \varphi_1(t, x)}{\partial t}  + \inf_{v \in V} \Bigl\{ c_1\bigl(t, x, v\bigr) + \mathcal{L}_{t}^{v} \varphi_1(t, x) \\
  + g_1\bigl(t, D_x \varphi_1(t, x) \cdot \sigma(t, x, (v,S(v)))\bigr)\Bigr\} = 0
\end{array}\right\}  \label{Eq3.5}
\end{eqnarray}
and 
\begin{eqnarray}
\left.\begin{array}{r}
 \dfrac{\partial \varphi_2(t, x)}{\partial t}  + \inf_{w \in W} \Bigl\{ c_2\bigl(t, x, w\bigr)  + \hat{\mathcal{L}}_{t}^{w} \varphi_2(t, x)  \\
  + g_2\bigl(t, D_x \varphi_2(t, x) \cdot \sigma(t, x, (\hat{v},w))\bigr)\Bigr\} = 0  \\
   \text{where} \,\, \hat{v} \,\, \text{is assumed to be given} 
\end{array}\right\}  \label{Eq3.6}
\end{eqnarray}
with, respectively, the following boundary conditions
\begin{align}
\varphi_1(T, x) = \Psi_1(T, x) \,\, \text{and} \,\, \varphi_2(T, x) = \Psi_2(T, x), \,\, x \in \mathbb{R}^d. \label{Eq3.7}
\end{align}

Note that the above equations in \eqref{Eq3.5} and \eqref{Eq3.6} together with \eqref{Eq3.7}, are associated with the risk-averse control problem w.r.t. the {\it leader} and {\it follower}, restricted to $\Sigma_{[t,T]}$, with cost functionals in \eqref{Eq2.17} and \eqref{Eq2.18}. Moreover, they represent generalized HJB equations with additional terms $g_1$ and $g_2$, respectively. Note that the problem of FBSDEs (cf. equations~\eqref{Eq2.1}, \eqref{Eq2.8} and \eqref{Eq2.9}) and the solvability of the corresponding HJB partial differential equations have been well studied in literature (e.g., see \cite{Ant93}, \cite{HUP95}, \cite{LIW14}, \cite{MaPY94}, \cite{ParT99}, \cite{Pen91} and \cite{Pen92}).

Next, we recall the definitions of viscosity solutions for \eqref{Eq3.5} and \eqref{Eq3.6} together with \eqref{Eq3.7} (e.g., see \cite{CraIL92}, \cite{FleS06} or \cite{Kry08} for additional discussions on the notion of viscosity solutions).

\begin{definition} \label{Df3}
The functions $\varphi_i \colon [0, T] \times \mathbb{R}^d$, with $i=1, 2$, are viscosity solutions for \eqref{Eq3.5} and \eqref{Eq3.6} together with the boundary conditions in \eqref{Eq3.7}, if the following conditions hold
\begin{enumerate} [(i)]
\item for every $\psi_i \in C_b^{1,2}([0, T], \times \mathbb{R}^d)$, with $i=1, 2$, such that $\psi_i \ge \varphi_i$ on $[0, T] \times \mathbb{R}^d$,
\begin{align}
\sup_{(t,x)} \bigl\{\varphi_i(t,x) - \psi_i(t,x) \bigr\} = 0, \label{Eq3.8}
\end{align}
and for $(t_{0},x_{0}) \in  [0, T] \times \mathbb{R}^d$ such that  $\psi_i(t_{0}, x_{0})=\varphi_i(t_{0}, x_{0})$ (i.e., a local maximum at $(t_{0},x_{0})$), then we have
\begin{align}
 &\dfrac{\partial \psi_1(t_{0},x_{0})}{\partial t}  + \inf_{v \in V} \Bigl\{ c_1\bigl(t_{0},x_{0}, v\bigr)  + \mathcal{L}_{t}^{v} \psi_1(t_{0},x_{0}) \notag \\
 &\quad  + g_1\bigl(t_{0}, D_x \psi_1(t_{0},x_{0}) \cdot \sigma(t_{0},x_{0},(v,S(v)))\bigr)\Bigr\} \ge 0\label{Eq3.9}
\end{align}
and 
\begin{align}
 &\dfrac{\partial \psi_2(t_{0},x_{0})}{\partial t}  + \inf_{w \in W} \Bigl\{ c_2\bigl(t_{0},x_{0}, w\bigr)  + \hat{\mathcal{L}}_{t}^{w} \psi_2(t_{0},x_{0}) \notag \\
 &\quad  + g_2\bigl(t_{0}, D_x \psi_2(t_{0},x_{0}) \cdot \sigma(t_{0},x_{0}, (v,w))\bigr)\Bigr\} \ge 0 \label{Eq3.10}
\end{align}
\item for every $\psi_i \in C_b^{1,2}([0, T], \times \mathbb{R}^d)$, with $i=1, 2$, such that $\psi_i \le \varphi_i$ on $[0, T] \times \mathbb{R}^d$,
\begin{align}
\inf_{(t,x)} \bigl\{ \varphi_i(t,x) - \psi_i(t,x) \bigr\} = 0,  \label{Eq3.11}
\end{align}
and for $(t_{0},x_{0}) \in  [0, T] \times \mathbb{R}^d$ such that  $\psi_i(t_{0}, x_{0})=\varphi_i(t_{0}, x_{0})$ (i.e., a local minimum at $(t_{0},x_{0})$), then we have
\begin{align}
 & \dfrac{\partial \psi_1(t_{0},x_{0})}{\partial t}  + \inf_{v \in V} \Bigl\{ c_1\bigl(t_{0},x_{0}, v\bigr) + \mathcal{L}_{t}^{v} \psi_1(t_{0},x_{0}) \notag \\
 &\quad  + g_1\bigl(t_{0}, D_x \psi_1(t_{0},x_{0}) \cdot \sigma(t_{0},x_{0}, (v,S(v)))\bigr)\Bigr\} \le 0 \label{Eq3.12}
\end{align}
and 
\begin{align}
& \dfrac{\partial \psi_2(t_{0},x_{0})}{\partial t}  + \inf_{w \in W} \Bigl\{ c_2\bigl(t_{0},x_{0}, w\bigr) + \hat{\mathcal{L}}_{t}^{w} \psi_2(t_{0},x_{0}) \notag \\
 & \quad  + g_2\bigl(t_{0}, D_x \psi_2(t_{0},x_{0}) \cdot \sigma(t_{0},x_{0}, (v,w))\bigr)\Bigr\} \le 0. \label{Eq3.13}
\end{align}
\end{enumerate}
\end{definition}

\subsection{On the risk-averse optimality condition for the follower} \label{S3.1}
{Note that, for a given $\hat{v}_{\cdot} \in \mathcal{V}_{[t,T]}$ (i.e., a risk-averse strategy for the} {\it leader's} {with restriction to $\Sigma_{[t, T]}$), if the risk-averse strategy for the} {\it follower} {is an optimal solution to} \eqref{Eq2.7} (cf. equations~\eqref{Eq2.4} and \eqref{Eq2.16}), then such a solution is characterized by the following propositions (i.e., Propositions~\ref{P3}, \ref{P4} and \ref{P5}). 

\begin{proposition} \label{P3}
Suppose that Assumption~\ref{AS1} together with \eqref{Eq1.7} and \eqref{Eq1.8} hold. Let $\hat{v}_{\cdot} \in \mathcal{V}_{[t,T]}$ be given, then the risk-value function w.r.t. the {\it follower} is given by
\begin{align}
 V_2^{w}\bigl(t, x\bigr) = \inf_{w_{\cdot} \in \mathcal{W}_{[t,r]} \bigl \vert \Sigma_{[t, T]}}  \rho_{t, r}^{g_2} \Bigl[\int_t^r c_2\bigl(s, X_s^{t,x;u}, w_s\bigr) ds + V_2^{w}\bigl(r, X_r^{t,x;u}\bigr) \Bigr] \label{Eq3.14}
\end{align}
 for any $(t,x) \in [0, T] \times \mathbb{R}^d$ and $r \in [t, T]$, with $u=(\hat{v}, w)$.
\end{proposition}

\begin{IEEEproof}
Notice that $\hat{v}_{\cdot} \in \mathcal{V}_{[t,T]}$ is given. Then, for any $\epsilon > 0$, there exists $\tilde{w}_{\cdot} \in \mathcal{W}_{[t,T]}$ such that $V_2^{w}\bigl(t, x\bigr)  + \epsilon \ge V_2^{\tilde{w}}\bigl(t, x\bigr)$. Further, if we applying the properties of time-consistency and translation (cf. Property~\ref{Pr1}) to $V_2^{\tilde{w}}\bigl(t, x\bigr)$, then we have
\begin{align}
V_2^{w}\bigl(t, x\bigr) + \epsilon \ge V_2^{\tilde{w}}\bigl(t, x\bigr)&= \rho_{t, r}^{g_2} \Bigl[\rho_{r, T}^{g_2} \Bigl[\int_t^T c_2\bigl(s, X_s^{t,x;\tilde{u}}, \tilde{w}_s\bigr) ds + \Psi_2(X_T^{t,x;\tilde{u}}) \Bigr] \Bigr] \notag \\
                                                    &= \rho_{t, r}^{g_2} \Bigl[ \int_t^r c_2\bigl(s, X_s^{t,x;\tilde{u}}, \tilde{w}_s\bigr) ds \notag \\
                                                    & \quad \quad + \rho_{r, T}^{g_2} \Bigl[ \int_r^T c_2\bigl(s, X_s^{t,x;\tilde{u}}, \tilde{w}_s\bigr) ds + \Psi_2(X_T^{t,x;\tilde{u}}) \Bigr] \Bigr], 
                                                    \label{EqP3.1}
\end{align}
where $\tilde{u}_{\cdot}=(\hat{v}_{\cdot}, \tilde{w}_{\cdot})$ is restricted to $\Sigma_{[t, T]}$. Moreover, if we apply Proposition~\ref{P2}, then we have  
\begin{align}
V_2^{w}\bigl(t, x\bigr) + \epsilon &\ge \rho_{t, r}^{g_2} \Bigl[ \int_t^r c_2\bigl(s, X_s^{t,x;\tilde{u}}, \tilde{w}_s\bigr) ds + V_2^{\tilde{w}}\bigl(r, X_r^{t,x;\tilde{u}}\bigr)\Bigr] \notag \\
                                                   &\ge \rho_{t, r}^{g_2} \Bigl[ \int_t^r c_2\bigl(s, X_s^{t,x;\tilde{u}}, \tilde{w}_s\bigr) ds + V_2^{w}\bigl(r, X_r^{t,x;\tilde{u}}\bigr)\Bigr] \notag \\
                                                   &\ge \inf_{w_{\cdot} \in \mathcal{W}_{[t,r]} \bigl \vert \Sigma_{[t, T]}} \rho_{t, r}^{g_2} \Bigl[ \int_t^r c_2\bigl(s, X_s^{t,x;\tilde{u}}, \tilde{w}_s\bigr) ds + V_2^{w}\bigl(r, X_r^{t,x;u}\bigr)\Bigr]. \label{EqP3.2}
\end{align}
Since $\epsilon$ is arbitrary, we obtain \eqref{Eq3.14}. On the other hand, to show the reverse inequality $``\le"$, let $\tilde{w}_{\cdot}$ (which is restricted to $\Sigma_{[t, T]}$) be an $\epsilon$-optimal solution, for some fixed $\epsilon > 0$, to the the problem on the right-hand side of \eqref{Eq3.14}.That is,
\begin{align}
\inf_{w_{\cdot} \in \mathcal{W}_{[t,r]} \bigl \vert \Sigma_{[t, T]}} \rho_{t, r}^{g_2} \Bigl[ \int_t^r c_2\bigl(s, X_s^{t,x;\tilde{u}}, \tilde{w}_s\bigr) ds + V_2^{w}\bigl(r, X_r^{t,x;u}\bigr)\Bigr] \notag \\
+ \epsilon  \quad \ge \rho_{t, r}^{g_2} \Bigl[ \int_t^r c_2\bigl(s, X_s^{t,x;\tilde{u}}, \tilde{w}_s\bigr) ds + V_2^{w}\bigl(r, X_r^{t,x;\tilde{u}}\bigr)\Bigr].  \label{EqP3.3}
\end{align}
Then, for every $y \in \mathbb{R}^d$, let $\tilde{w}_{\cdot}(y) \in \mathcal{W}_{[t,T]}$ be such that $V_2^{w}\bigl(r, y\bigr) + \epsilon \ge V_2^{\tilde{w}(y)}\bigl(t, x\bigr)$ and restricted to $\Sigma_{[t, T]}$. Due to the measurable selection theorem (e.g., see \cite[Theorem~I.7.7,~pp.~161]{Wag72}), we may assume that the function $y \rightarrow \tilde{w}(y)$ is Borel measurable. Further, suppose that a control function $w_{\cdot}^0$ is defined as follow
\begin{eqnarray}
w_s^0 =\left\{\begin{array}{l l}
 \bar{w}_s,  & s \in [t, r)\\
  \tilde{w}_s(X_s^{t,x;\bar{u}}), & s \in [r, T].
\end{array}\right.  \label{EqP3.4}
\end{eqnarray}
Note that, from the above definition, $w_{\cdot}^0$ is restricted to $\Sigma_{[t, T]}$. Then, using the properties of the monotonicity, translation and time-consistency, we obtain the following
\begin{align}
 \rho_{t, r}^{g_2} &\Bigl[ \int_t^r c_2\bigl(s, X_s^{t,x;\bar{u}}, \bar{w}_s\bigr) ds + V_2^{\bar{w}}\bigl(r, X_r^{t,x;\bar{u}}\bigr)\Bigr] \notag \\
                                                   &\ge \rho_{t, r}^{g_2} \Bigl[ \int_t^r c_2\bigl(s, X_s^{t,x;\bar{u}}, \bar{w}_s\bigr) ds + V_2^{\tilde{w}_s(X_s^{t,x;\bar{u}})}\bigl(r, X_r^{t,x;\bar{u}}\bigr) -\epsilon \Bigr], \notag \\
                                                   & \quad \quad  \text{with} \,\, \bar{u} =(\hat{v}, \bar{w}) \notag \\
                                                   &\ge \rho_{t, T}^{g_2} \Bigl[ \int_t^T c_2\bigl(s, X_s^{t,x;u^0}, \bar{w}_s^0\bigr) ds + \Psi_2\bigl(X_T^{t,x;u^0}\bigr)\Bigr] -\epsilon, \notag \\
                                                   & \quad \quad  \text{with} \,\, u^0 =(\hat{v}, w^0)  \notag \\
                                                   &= V_2^{w^{0}} \bigl(t, x\bigr) - \epsilon. \label{EqP3.5}
\end{align}
If we further combine the inequalities from \eqref{EqP3.3} and \eqref{EqP3.5}, then we have
\begin{align}
\inf_{w_{\cdot} \in \mathcal{W}_{[t,r]} \bigl \vert \Sigma_{[t, T]}} \rho_{t, r}^{g_2} \Bigl[  \int_t^r c_2\bigl(s, X_s^{t,x;u}, w_s\bigr) ds + V_2^{w}\bigl(r, X_r^{t,x;u}\bigr)\Bigr]  + \epsilon & \ge V_2^{w^0}\bigl(t, x\bigr) - \epsilon \notag \\
  & \ge V_2^{w}\bigl(t, x\bigr) - \epsilon.  \label{EqP3.6}
\end{align}
Note that, since $\epsilon$ is arbitrary, we obtain \eqref{Eq3.14}. This completes the proof of Proposition~\ref{P3}. 
\end{IEEEproof}

Then, we have the following results (i.e., Propositions~\ref{P4} and \ref{P5}) {that are useful to characterize the mapping $S$ in} \eqref{Eq2.12}.
\begin{proposition}\label{P4}
Suppose that Assumption~\ref{AS1} holds and let $W$ be a compact set in $\mathbb{R}^{d}$. Let $\hat{v}_{\cdot} \in \mathcal{V}_{[t,T]}$ be given, then the risk-value function $V_2^{w}\bigl(\cdot, \cdot\bigr)$ is the viscosity solution of \eqref{Eq3.6} with boundary condition $\Psi_2(T, x)$ for $x \in \mathbb{R}^d$ and with $u=(\hat{v}, w)$.
\end{proposition}

\begin{IEEEproof}
Suppose that $\varphi_2 \in C_b^{1,2}([0, T] \times \mathbb{R}^d)$ and assume that $\varphi_2 \ge V_2^{w}$ on $[0, T] \times \mathbb{R}^d$ and $\max_{(t,x)} \bigl[V_2^{w}(t,x) - \varphi_2(t,x)\bigr] = 0$. We consider a point $(t_{0},x_{0}) \in [0, T] \times \mathbb{R}^d$ so that $\varphi_2(t_{0},x_{0})= V_2^{w}(t_{0},x_{0})$ (i.e., a local maximum at $(t_{0},x_{0})$). Further, for a small $\delta t > 0$, we consider a constant control $w_s=\alpha$ for $s \in [t_{0},t_{0} + \delta t]$. Then, from \eqref{Eq3.14}, we have
\begin{align}
 \varphi_2(t_{0},x_{0}) &= V_2^{w}(t_{0},x_{0}) \notag \\
                                        & \le \rho_{t_{0},t_{0} + \delta t}^{g_2} \Bigl[\int_{t_{0}}^{t_{0} + \delta t} c_2\bigl(s, X_s^{t_{0},x_{0};u}, \alpha \bigr) ds + V_2^{w}(t_{0} + \delta t, X_{t_{0} + \delta t}^{t_{0},x_{0};u})\Bigr]  \notag \\
                                       & \le \rho_{t_{0},t_{0} + \delta t}^{g_2} \Bigl[\int_{t_{0}}^{t_{0} + \delta t} c_2\bigl(s, X_s^{t_{0},x_{0};u}, \alpha \bigr) ds + \varphi_2(t_{0} + \delta t, X_{t_{0} + \delta t}^{t_{0},x_{0};u})\Bigr], \quad  \text{with} \,\, u=(\hat{v}, \alpha). \label{EqP4.1}
\end{align}
Using the translation property of $\rho_{t_{0},t_{0} + \delta t}[\,\cdot\,]$, we obtain the following inequality
\begin{align}
\rho_{t_{0},t_{0} + \delta t}^{g_2} \Bigl[\int_{t_{0}}^{t_{0} + \delta t} c_2\bigl(s, X_s^{t_{0},x_{0};u}, \alpha \bigr) ds + \varphi_2(t_{0} + \delta t, X_{t_{0} + \delta t}^{t_{0},x_{0};u}) - \varphi_2(t_{0},x_{0})\Bigr] \ge 0. \label{EqP4.2}
\end{align}
Notice that $\varphi_2 \in C_b^{1,2}([0, T] \times \mathbb{R}^d)$, then, using the It\^{o} formula, we can evaluate the difference between $\varphi_2(t_{0} + \delta t, X_{t_{0} + \delta t}^{t_{0},x_{0};u})$ and $\varphi_2(t_{0},x_{0})$ as follow
\begin{align}
&\varphi_2(t_{0} + \delta t, X_{t_{0} + \delta t}^{t_{0},x_{0};u}) - \varphi_2(t_{0},x_{0}) \notag \\
&= \int_{t_{0}}^{t_{0} + \delta t} \Bigl[\dfrac{\partial}{\partial t} \varphi_2(s, X_{s}^{t_{0},x_{0};u}) + \hat{\mathcal{L}}_{t}^{\alpha} \varphi_2(s, X_{s}^{t_{0},x_{0};u}) \Bigr] d s \notag  \\
  & \quad + \int_{t_{0}}^{t_{0} + \delta t} D_x \varphi_2(s, X_{s}^{t_{0},x_{0};u}) \cdot \sigma(s, X_{s}^{t_{0},x_{0};u}, (\hat{v}_s,\alpha)) d B_s. \label{EqP4.3}
\end{align}
Moreover, if we substitute the above equation into \eqref{EqP4.2}, then we obtain
\begin{align}
&\rho_{t_{0},t_{0} + \delta t}^{g_2} \Bigl[ \int_{t_{0}}^{t_{0} + \delta t} \Bigl[c_2\bigl(s, X_s^{t_{0},x_{0};u}, \alpha \bigr) + \dfrac{\partial}{\partial t} \varphi_2(s, X_{s}^{t_{0},x_{0};u}) + \hat{\mathcal{L}}_{t}^{\alpha} \varphi_2(s, X_{s}^{t_{0},x_{0};u}) \Bigr] d s\notag \\
& \quad \quad + \int_{t_{0}}^{t_{0} + \delta t} D_x \varphi_2(s, X_{s}^{t_{0},x_{0};u}) \cdot \sigma(s, X_{s}^{t_{0},x_{0};u}, (\hat{v}_s,\alpha)) d B_s \Bigr]  \ge 0, \label{EqP4.4}
\end{align}
which amounts to solving the following BSDE
\begin{align}
& Y_{t_{0}}^{2;t_{0},x_{0};u} = \int_{t_{0}}^{t_{0} + \delta t} \Bigl[c_2\bigl(s, X_s^{t_{0},x_{0};u}, \alpha \bigr) + \dfrac{\partial}{\partial t} \varphi_2(s, X_{s}^{t_{0},x_{0};u}) + \hat{\mathcal{L}}_{t}^{\alpha} \varphi_2(s, X_{s}^{t_{0},x_{0};u}) \Bigr] d s \notag \\
& \quad \quad + \int_{t_{0}}^{t_{0} + \delta t} D_x \varphi_2(s, X_{s}^{t_{0},x_{0};u}) \cdot \sigma(s, X_{s}^{t_{0},x_{0};u}, (\hat{v}_s,\alpha)) d B_s \notag \\
&\quad \quad \quad + \int_{t_{0}}^{t_{0} + \delta t}g_2\bigl(s, Z_s^{2;t_{0},x_{0};u}\bigr) ds - \int_{t_{0}}^{t_{0} + \delta t} Z_s^{2;t_{0},x_{0};u} dB_s.  \label{EqP4.5}
\end{align}
From Lemma~\ref{L1}, the above BSDE admits unique solutions, i.e.,
\begin{align*}
Z_s^{2;t_{0},x_{0};u} = D_x \varphi_2(s, X_{s}^{t_{0},x_{0};u}) \cdot \sigma(s, X_{s}^{t_{0},x_{0};u}, (\hat{v}_s,\alpha)), \,\, t_{0} \le s \le t_{0} + \delta t
\end{align*}
and
\begin{align*}
Y_{t_{0}}^{2;t_{0},x_{0};u}  = \int_{t_{0}}^{t_{0} + \delta t} \Bigl[c_2\bigl(s, X_s^{t_{0},x_{0};u}, \alpha \bigr) &+ \dfrac{\partial}{\partial t} \varphi_2(s, X_{s}^{t_{0},x_{0};u}) + \hat{\mathcal{L}}_{t}^{\alpha} \varphi_2(s, X_{s}^{t_{0},x_{0};u}) \notag \\
&  \quad    + g_2\bigl(s, D_x \varphi_2(s, X_{s}^{t_{0},x_{0};u}) \cdot \sigma(s, X_{s}^{t_{0},x_{0};u}, (\hat{v}_s,\alpha))\bigr) \Bigr] d s.
\end{align*}
Further, if we substitute the above results in \eqref{EqP4.4}, we obtain
\begin{align}
 \int_{t_{0}}^{t_{0} + \delta t} \Bigl[c_2\bigl(s, X_s^{t_{0},x_{0};u}, \alpha \bigr) &+ \dfrac{\partial}{\partial t} \varphi_2(s, X_{s}^{t_{0},x_{0};u}) + \hat{\mathcal{L}}_{t}^{\alpha} \varphi_2(s, X_{s}^{t_{0},x_{0};u})  \notag \\
& + g_2\bigl(s, D_x \varphi_2(s, X_{s}^{t_{0},x_{0};u}) \cdot \sigma(s, X_{s}^{t_{0},x_{0};u}, (\hat{v}_s,\alpha))\bigr)  \Bigr] d s \ge 0.  \label{EqP4.6}
\end{align}
Then, dividing the above equation by $\delta t$ and letting $\delta t \rightarrow 0$, we obtain
\begin{align*}
 c_2\bigl(t_{0},x_{0}, \alpha \bigr) & + \dfrac{\partial}{\partial t} \varphi_2(t_{0},x_{0}) + \hat{\mathcal{L}}_{t}^{\alpha} \varphi_2(t_{0},x_{0})+ g_2\bigl(t_{0}, D_x \varphi_2(t_{0},x_{0}) \cdot \sigma(t_{0},x_{0}, (\hat{v}_{t_{0}},\alpha))\bigr) \ge 0. 
\end{align*}
Note that, since $\alpha \in W$ is arbitrary, we can rewrite the above condition as follow
\begin{align}
& \dfrac{\partial}{\partial t_{0}} \varphi_2(t_{0},x_{0}) + \min_{\alpha \in W} \Bigl\{c_2\bigl(t_{0},x_{0}, \alpha \bigr) + \hat{\mathcal{L}}_{t}^{\alpha} \varphi_2(t_{0},x_{0})  + g_2\bigl(t_{0}, D_x \varphi_2(t_{0},x_{0}) \cdot \sigma(t_{0},x_{0}, (\hat{v}_{t_{0}},\alpha))\bigr) \Bigr\} \ge 0, \label{EqP4.7} 
\end{align}
which attains its minimum in $W$ (which is a compact set in $\mathbb{R}^d)$. Thus, $V_2^{w}(\cdot,\cdot)$ is a viscosity subsolution of \eqref{Eq3.15}, with boundary condition $\varphi_2(T, x)=\Psi_2(T, x)$.

On the other hand, suppose that $\varphi_2 \in C_b^{1,2}([0, T] \times \mathbb{R}^d)$ and assume that $\varphi_2 \le V_2^{w}$ on $[0, T] \times \mathbb{R}^d$ and $\min_{(t,x)} \bigl[V_2^{w}(t,x) - \varphi_2(t,x)\bigr] =0$. Then, we consider a point $(t_{0},x_{0}) \in [0, T] \times \mathbb{R}^d$ so that $\varphi_2(t_{0},x_{0})= V_2^{w}(t_{0},x_{0})$ (i.e., a local minimum at $(t_{0},x_{0})$). Further, for a small $\delta t > 0$, Let $\tilde{w}_s$, which is restricted to $\Sigma_{[t_{0},t_{0} + \delta t]}$, be an $\epsilon \delta t$-optimal control for \eqref{Eq3.14} at $(t_{0},x_{0})$. Then, proceeding in this way as \eqref{EqP4.6}, we obtain the following 
\begin{align}
\int_{t_{0}}^{t_{0} + \delta t}  \Bigl[c_2\bigl(s, X_s^{t_{0},x_{0};u}, \tilde{w}_s \bigr)& + \dfrac{\partial}{\partial t} \varphi_2(s, X_{s}^{t_{0},x_{0};u}) + \hat{\mathcal{L}}_{t}^{\tilde{w}} \varphi_2(s, X_{s}^{t_{0},x_{0};u}) \notag \\
&\,\, + g_2\bigl(s, D_x \varphi_2(s, X_{s}^{t_{0},x_{0};u}) \cdot \sigma(s, X_{s}^{t_{0},x_{0};u}, (\hat{v}_s,\tilde{w}_s))\bigr) \Bigr] d s \le \epsilon \delta t, \notag \\
& \quad \quad  \text{with} \,\,   u =(\hat{v}_s, \tilde{w}_s).  \label{EqP4.8}
\end{align}
As a result of this, we also obtain the following
\begin{align}
\int_{t_{0}}^{t_{0} + \delta t}  \min_{\alpha \in W} \Bigl\{c_2\bigl(s, X_s^{t_{0},x_{0};u}, \alpha \bigr) &+ \dfrac{\partial}{\partial t} \varphi_2(s, X_{s}^{t_{0},x_{0};u}) + \hat{\mathcal{L}}_{t}^{\alpha} \varphi_2(s, X_{s}^{t_{0},x_{0};u}) \notag \\
& \quad  + g_2\bigl(s, D_x \varphi_2(s, X_{s}^{t_{0},x_{0};u}) \cdot \sigma(s, X_{s}^{t_{0},x_{0};u}, (\hat{v}_s,\alpha))\bigr) \Bigr\} d s  \le \epsilon \delta t.  \label{EqP4.9}
\end{align}
Note that the mapping 
\begin{align*}
& (s, x, \alpha) \rightarrow  \Bigl[c_2\bigl(s, x, \alpha \bigr) + \dfrac{\partial}{\partial t} \varphi_2(t, x) + \hat{\mathcal{L}}_{t}^{\alpha} \varphi_2(t, x)   + g_2\bigl(t, D_x \varphi_2(t, x) \cdot \sigma(t, x, (\hat{v}_t,\alpha))\bigr) \Bigr]
\end{align*}
is continuous and, since $W$ is compact, then $s \rightarrow X_s^{t_{0},x_{0};u}$ is also continuous. As a result, the expression under the integral in \eqref{EqP4.9} is continuous. Further, if we divide both sides of \eqref{EqP4.9} by $\delta t$ and letting $\delta t \rightarrow 0$, then we obtain the following
\begin{align}
& \dfrac{\partial}{\partial t_{0}} \varphi_2(t_{0},x_{0})  + \min_{\alpha \in W} \Bigl\{c_2\bigl(t_{0},x_{0}, \alpha \bigr) + \hat{\mathcal{L}}_{t}^{\alpha} \varphi_2(t_{0},x_{0})  + g_2\bigl(t_{0}, D_x \varphi_2(t_{0},x_{0}) \cdot \sigma(t_{0},x_{0}, (\hat{v}_{t_{0}},\alpha))\bigr) \Bigr\} \le \epsilon. \label{EqP4.10} 
\end{align}
Notice that, since $\epsilon$ is arbitrary, we conclude that $V_2^{w}(\cdot,\cdot)$ is a viscosity supersolution of \eqref{Eq3.13}, with boundary condition $\varphi_2(T,x)=\Psi_2(T, x)$. This completes the proof of Proposition~\ref{P4}.
\end{IEEEproof}

\begin{remark} \label{R5}
Note that if \,$V_2^{w} \in C_b^{1,2}([0, T] \times \mathbb{R}^d)$, then such a solution also satisfies \eqref{Eq3.13} with boundary condition $V_2^{w}(T,x)=\Psi_2(T, x)$. Furthermore, using the {\it verification theorem}, one can also identify $V_2^{w}$ as the optimal value function.
\end{remark}

\begin{proposition} \label{P5}
Suppose that Proposition~\ref{P4} holds and let $\varphi_2 \in C_b^{1,2}([0, T] \times \mathbb{R}^d)$ satisfy \eqref{Eq3.6} with $\varphi_2\bigl(T, x\bigr)=\Psi_2(T, x)$ for $x \in \mathbb{R}^d$. Then, $\varphi_2\bigl(t, x\bigr) \le V_2^{w}\bigl(t, x\bigr)$ for any control $w_{\cdot} \in \mathcal{W}_{[t,T]}$ with restriction to $\Sigma_{[t, T]}$ and for all $(t,x) \in [0, T] \times \mathbb{R}^d$. Furthermore, if an admissible control process $\hat{w}_{\cdot} \in \mathcal{W}_{[t,T]}$ exists, for almost all $(s, \Omega) \in [0, T] \times \Omega$, together with the corresponding solution $X_s^{t,x;\hat{u}}$, with $\hat{u}_s=(\hat{v}_s, \hat{w}_s)$, and satisfies
\begin{align}
 \hat{w}_s \in & \arginf_{w_{\cdot} \in \mathcal{W}_{[t,T]} \bigl \vert \Sigma_{[t, T]}} \Bigl\{c_2\bigl(s, X_s^{t,x;u}, w_s\bigr) + \hat{\mathcal{L}}_{s}^{w} \varphi_2\bigl(s, X_s^{t,x;u}\bigr) \notag \\
 & \underbrace{\hspace{0.25in}  + g_2\bigl(s, D_x\varphi_2 \bigl(s, X_s^{t,x;u}\bigr) \cdot \sigma\bigl(s, X_s^{t,x;u}, (\hat{v}_s, w_s\bigr)\bigr)\bigr)\Bigr\}}_{\triangleq S(\hat{v}) \,\, \text{with} \,\, S \colon \mathcal{V}_{[t,T]} \rightrightarrows \mathcal{W}_{[t,T]}} \label{Eq3.15}
 \end{align}
Then, $\varphi_2\bigl(t, x\bigr) = V_2^{\hat{w}}\bigl(t, x\bigr)$ for all $(t,x) \in [0, T] \times \mathbb{R}^d$.
\end{proposition}

\begin{IEEEproof}
Assume that $(t, x) \in [0, T] \times \mathbb{R}^d$ is fixed. For any $w_{\cdot} \in \mathcal{W}_{[t,T]}$, restricted to $\Sigma_{[t, T]}$, we consider a process $\kappa \bigl(s, X_{s}^{t,x;u}\bigr)$, with $u=(\hat{v}, w)$, for $s \in [t, T]$. Then, using It\^{o} integral formula, we can evaluate the difference between $\kappa \bigl(T, X_{T}^{t,x;u}\bigr)$ and $\kappa \bigl(t, x\bigr)$ as follow\footnote{Notice that $\kappa \bigl(t, x\bigr) \in C_b^{1,2}([0, T] \times \mathbb{R}^d).$}
\begin{align}
\kappa \bigl(T, X_{T}^{t,x;u}\bigr) - \kappa \bigl(t, x\bigr) &= \int_{t}^T \Bigl[\dfrac{\partial}{\partial t}\kappa\bigl(s, X_{s}^{t,x;u}\bigr) + \hat{\mathcal{L}}_{t}^{w} \kappa \bigl(s, X_{s}^{t,x;u}\bigr) \Bigr] d s \notag  \\
  & \quad + \int_{t}^T D_x \kappa \bigl(s, X_{s}^{t,x;u}\bigr) \cdot \sigma(s, X_{s}^{t,x;u}, (\hat{v}_s,w_s)) d B_s. \label{EqP5.1}
\end{align}
Using \eqref{Eq3.6}, we further obtain the following
\begin{align}
&  \dfrac{\partial}{\partial t} \kappa\bigl(s, X_{s}^{t,x;u}\bigr) + \hat{\mathcal{L}}_{t}^{w} \kappa\bigl(s, X_{s}^{t,x;u}\bigr) + c_2\bigl(s, X_{s}^{t,x;u}, w_s \bigr) \notag \\
& \quad\quad + g_2\bigl(s, D_x \varphi_2(s, X_{s}^{t,x;u}) \cdot \sigma(s, X_{s}^{t,x;u}, (\hat{v}_s,w_s))\bigr) \ge 0. \label{EqP5.2}
\end{align}
Further, if we combine \eqref{EqP5.1} and \eqref{EqP5.2}, then we obtain
\begin{align}
& \kappa \bigl(t, x\bigr) \le \Psi_2\bigl(T, X_{T}^{t,x;u}\bigr) + \int_{t}^T \Bigl\{c_2\bigl(s, X_{s}^{t,x;u}, w_s \bigr) ds \notag \\
                                    & \quad + g_2\bigl(s, D_x \kappa(s, X_{s}^{t,x;u}) \cdot \sigma(s, X_{s}^{t,x;u}, (\hat{v}_s,w_s))\bigr) \Bigr\} ds \notag \\
                                     & \quad -  \int_{t}^T D_x \kappa(s, X_{s}^{t,x;u}) \cdot \sigma(s, X_{s}^{t,x;u}, (\hat{v}_s,w_s)) d B_s. \label{EqP5.3}
\end{align}
Define $Z_{s}^{2;t,x;u} = D_x \kappa(s, X_{s}^{t,x;u}) \cdot \sigma(s, X_{s}^{t,x;u}, (\hat{v}_s,w_s))$, for $s \in [t, T]$, then $\kappa \bigl(t, x\bigr) \le Y_{t}^{2, t,x;u}$ follows, where $(Y_{\cdot}^{2;t,x;u}, Z_{\cdot}^{2;t,x;u})$ is a solution to BSDE in \eqref{Eq2.9}. As a result of this, we have
\begin{align*}
 \kappa \bigl(t, x\bigr) \le V_2^{w}\bigl(t, x\bigr).
\end{align*}
Moreover, if there exists at least one $\hat{w}$ satisfying \eqref{Eq3.15}, then, for $w=\hat{w}$, the inequality in \eqref{EqP5.3} becomes an equality (i.e., $\kappa(t,x)=V_2^{\hat{w}}\bigl(t, x\bigr)$). Note that the corresponding pathwise solution $X_{s}^{t,x;\hat{u}}$, with $\hat{u}=(\hat{v}, \hat{w})$ and $\hat{w}=S(\hat{v})$, is progressively measurable, since the control process $\hat{w}_{\cdot}  \in \mathcal{W}_{[t,T]}$ is restricted to $\Sigma_{[t, T]}$. This completes the proof of Proposition~\ref{P5}. 
\end{IEEEproof}

\subsection{On the risk-averse optimality condition for the leader} \label{S3.2}
In this subsection, we provide an optimality condition on the risk-averseness of the {\it leader} in \eqref{Eq2.13} (cf. equations~\eqref{Eq2.15} and \eqref{Eq2.2}), when the risk-averse strategy for the {\it follower} satisfies the optimality condition of Proposition~\ref{P5}.

\begin{proposition} \label{P6}
Suppose that Proposition~\ref{P5} holds. Then, for any $(t,x) \in [0, T] \times \mathbb{R}^d$ and $r \in [t, T]$, the risk-value function w.r.t. the {\it leader} is given by
\begin{align}
 V_1^{v}\bigl(t, x\bigr) = \inf_{v_{\cdot} \in \mathcal{V}_{[t,r]}\bigl \vert \Sigma_{[t, T]}}  \rho_{t, r}^{g_1} \Bigl[\int_t^r c_1\bigl(s, X_s^{t,x;u}, v_s\bigr) ds + V_1^{v}\bigl(r, X_r^{t,x;u}\bigr) \Bigr] \label{Eq3.16}
\end{align}
Furthermore, if  \,$V$ is a compact set in $\mathbb{R}^{d}$, then $V_1^{v}\bigl(\cdot, \cdot\bigr)$ is the viscosity solution of \eqref{Eq3.5} with boundary condition $\Psi_1(T, x)$ for $x \in \mathbb{R}^d$.
\end{proposition}

\begin{IEEEproof}
Note that if Proposition~\ref{P5} holds, then, for any control process $v_{\cdot} \in \mathcal{V}_{[t,T]}$, there exists at least one $w_{\cdot} \in S(v_{\cdot}) \subset \mathcal{W}_{[t,T]}$ such that $S$ satisfies the definition of the mapping in \eqref{Eq3.15}. As a result, one can prove the above proposition {(w.r.t. the operator $\mathcal{L}_{t}^{v}$)} in the same way as Proposition~\ref{P4}. 
\end{IEEEproof}

\begin{proposition} \label{P7}
Suppose that Proposition~\ref{P6} holds and let $\varphi_1 \in C_b^{1,2}([0, T] \times \mathbb{R}^d)$ satisfy \eqref{Eq3.5} with boundary condition $\varphi_1\bigl(T, x\bigr)=\Psi_1(T, x)$ for $x \in \mathbb{R}^d$. Then, $\varphi_1\bigl(t, x\bigr) \le V_1^{v}\bigl(t, x\bigr)$ for any control process $v_{\cdot} \in \mathcal{V}_{[t,T]}$ with restriction to $\Sigma_{[t, T]}$ and for all $(t,x) \in [0, T] \times \mathbb{R}^d$. Furthermore, if an admissible control process $v_{\cdot}^{\ast} \in \mathcal{V}_{[t,T]}$ exists, for almost all $(s, \Omega) \in [t, T] \times \Omega$, together with the corresponding solution $X_s^{t,x;u*}$, with $u_s^{\ast}=(v_s^{\ast}, S(v_s^{\ast}))$, and satisfies
\begin{align}
& v_s^{\ast} \in  \arginf_{v_{\cdot} \in \mathcal{V}_{[t,T]} \bigl \vert \Sigma_{[t, T]}} \Bigl\{c_1\bigl(s, X_s^{t,x;u}, v_s\bigr) + \mathcal{L}_{s}^{v} \varphi_1\bigl(s, X_s^{t,x;u}\bigr) \notag \\
 &\quad  \quad + g_1\bigl(s, D_x\varphi_1 \bigl(s, X_s^{t,x;u}\bigr) \cdot \sigma\bigl(s, X_s^{t,x;u}, (v_s, S(v_s)\bigr)\bigr) \bigr) \Bigr\}. \label{Eq3.17}
 \end{align}
Then, $\varphi_1\bigl(t, x\bigr) = V_1^{v*}\bigl(t, x\bigr)$ for all $(t,x) \in [0, T] \times \mathbb{R}^d$.
\end{proposition}

\begin{IEEEproof}
The proof is similar to that of Proposition~\ref{P5}, except that we require a unique solution set $\bigl\{X_{\cdot}^{t,x;u}, (Y_{\cdot}^{1;t,x;u}, Z_{\cdot}^{1;t,x;u}), (Y_{\cdot}^{2;t,x;u}, Z_{\cdot}^{2;t,x;u})\bigr\}$ for the FBSDEs in \eqref{Eq2.1}, \eqref{Eq2.8} and \eqref{Eq2.9} on $\bigl(\Omega, \mathcal{F}, P, \mathcal{F}^t\bigr)$ for every initial condition $(t,x) \in [0,T] \times \mathbb{R}^d$, with $u_{s}^{\ast}=(v_{s}^{\ast}, S(v_{s}^{\ast}))$ and $s \in [t, T]$.
\end{IEEEproof}

Here, we remark that the above condition in \eqref{Eq3.17} requires the {\it follower} to respond optimally to the risk-averse strategy of the {\it leader}, where such a correspondence is implicitly embedded via the mapping $S$ (cf. equation~\eqref{Eq3.15}). {Furthermore, the optimal risk-averse strategies $v_{\cdot}^{\ast}$ and $S(v_{\cdot}^{\ast})$ are interpreted in the sense of viscosity for the corresponding risk-averse dynamic programming equations in} \eqref{Eq3.14} and \eqref{Eq3.16}. 

\section{Further remarks} \label{S4}
In this section, {we briefly comment on the implication of our results} in assessing the influence of the {\it leader}'s decision and its risk-averse satisfaction on the risk-averseness of the {\it follower} in relation to the direction of {\it leader-follower} information flow. Note that, for the optimization problems in \eqref{Eq2.15} and \eqref{Eq2.16}, if there exist: {(i) a nonanticipating measurable mapping $S \colon \mathcal{V}_{[0,T]} \rightrightarrows \mathcal{W}_{[0,T]}$ that satisfies}
\begin{align*}
\rho_{0, T}^{g_2}\bigr[\xi_{0,T}^2\bigl(v, S(v)\bigr)\bigl] \le \rho_{0, T}^{g_2}\bigr[\xi_{0,T}^2\bigl(v, w\bigr)\bigl],\, v_{\cdot} \in \mathcal{V}_{[0,T]}, \,\, \forall w_{\cdot} \in \mathcal{W}_{[0,T]}
\end{align*}
{and (ii) a risk-averse strategy for the {\it leader} $v_{\cdot}^{\ast} \in \mathcal{V}_{[0,T]}$ that satisfies}
\begin{align*}
\rho_{0, T}^{g_1}\bigr[\xi_{0,T}^1\bigl(v^{\ast}, S(v^{\ast})\bigr)\bigl] \le \rho_{0, T}^{g_1}\bigr[\xi_{0,T}^1\bigl(v, S(v)\bigr)\bigl],\, \forall v_{\cdot} \in \mathcal{V}_{[0,T]}.
\end{align*}
{Then, the pair $(v_{\cdot}^{\ast}, w_{\cdot}^{\ast}) \in \mathcal{V}_{[0,T]} \times\mathcal{W}_{[0,T]}$ is an optimal risk-averse strategy. Moreover, we observe that the statement} of Proposition~\ref{P5} is implicitly accounted in Proposition~\ref{P7}. That is, for any $s \in [t, T]$, the risk-averseness of the {\it follower}, with restriction to $\Sigma_{[t,T]}$, i.e.,
\begin{align*}
w_s^{\ast} \in & \arginf_{w_{\cdot} \in \mathcal{W}_{[t,T]}\bigl \vert \Sigma_{[t, T]}} \Bigl\{c_2\bigl(s, X_s^{t,x;u}, w_s\bigr) + \hat{\mathcal{L}}_{s}^{w} \varphi_2\bigl(s, X_s^{t,x;u}\bigr)  \\
 &\quad +g_2\bigl(s, D_x\varphi_2 \bigl(s, X_s^{t,x;u}\bigr) \cdot \sigma\bigl(s, X_s^{t,x;u}, (v_s, w_s)\bigr)\bigr) \Bigr\}
\end{align*}
is a subproblem in \eqref{Eq3.17} and its solution is also given by
\begin{align*}
w_s^{\ast} \in & \arginf_{w_{\cdot} \in S_b(v_{\cdot}^{\ast}) \bigl \vert \Sigma_{[t, T]}} \Bigl\{c_1\bigl(s, X_s^{t,x;u}, v_s^{\ast}\bigr) + \mathcal{L}_{s}^{(v^{\ast}, w)} \varphi_1\bigl(s, X_s^{t,x;u}\bigr) \\
 &\quad +g_1\bigl(s, D_x\varphi_1 \bigl(s, X_s^{t,x;u}\bigr) \cdot \sigma\bigl(s, X_s^{t,x;u}, (v_s^{\ast}, w_s)\bigr)\bigr) \Bigr\},
\end{align*}
provided that the set of all {\it best-responses} of the {\it follower} $S_b(v_{\cdot}^{\ast})$, i.e., 
\begin{align*}
S_b(v_{\cdot}^{\ast}) \subseteq \bigl\{w_{\cdot} \in \mathcal{W}_{[t,T]} \,\bigl\vert\, v_{\cdot}^{\ast} \in \mathcal{V}_{[t,T]} \,\, \text{and} \,\, \rho_{t, T}^{g_2} \bigl[\xi_{t,T}^2\bigl(v^{\ast}, w\bigr)\bigl]  \le \rho_{t, T}^{g_2} \bigl[\xi_{t,T}^2\bigl(v^{\ast}, \hat{w}\bigr)\bigl], \,\, \forall \hat{w}_{\cdot} \in \mathcal{W}_{[t,T]} \bigr\}
\end{align*}
is nonempty. On the other hand, the risk-averse strategy for the {\it leader} is given by
\begin{align*}
  v_s^{\ast} \in & \arginf_{v_{\cdot} \in \mathcal{V}_{[t,T]} \bigl \vert \Sigma_{[t, T]}} \Bigl\{c_1\bigl(s, X_s^{t,x;u}, v_s\bigr) + \mathcal{L}_{s}^{v} \varphi_1\bigl(s, X_s^{0,x;u}\bigr) \\
 &  \quad + g_1\bigl(s, D_x\varphi_1 \bigl(s, X_s^{t,x;u}\bigr) \cdot \sigma\bigl(s, X_s^{t,x;u}, (v_s, S(v_s)\bigr)\bigr) \bigr) \Bigr\},
\end{align*}
{where such a strategy is also implicitly conditioned by the} {\it follower's} {risk-averse strategy and its risk-averse satisfaction.} As a result of this, the {\it follower} is involved not only in minimizing his own accumulated risk-cost (in response to the risk-averse strategy of the {\it leader}) but also in minimizing that of the {\it leader's} accumulated risk-cost. Hence, such an inherent interaction, due to the {\it nature of the problem}, constitutes a constrained information flow between the {\it leader} and that of the {\it follower}, in which the {\it follower} is required to respond optimally, in the sense of {\it best-response correspondence}, to the risk-averse strategy of the {\it leader}.

\begin{remark} [{\it Comment on model uncertainty}]\label{R8}
Finally, it is worth remarking that the issue of risk-averseness under model uncertainty, when the {\it follower} is allowed to take into account alternative models that are statistically difficult to distinguish from \eqref{Eq2.1}, can be modeled as an ``{\it min-max}" optimization problem with uncertainty aversion. {For example, if we replace $B_t$ with $(\widehat{B}_t + \int_0^t h_s d s)$ in} \eqref{Eq2.1}, {where $h \in C_b^{1, 2} \bigl([0, T]; \mathbb{R}^{d}\bigr)$ is a measurable function and $\widehat{B}_t$ is a $d$-dimensional brownian motion.} Then, the process $h$ can be used as a device to transform $\mathbb{P}$ into another probability distribution $\mathbb{Q}$ which is mutually absolutely continuous w.r.t. $\mathbb{P}$ on $\bigl(\Omega, \mathcal{F}\bigr)$ (e.g., see Girsavov \cite{Gir60} for additional discussions on transforming stochastic processes). Moreover, if we specify the model uncertainty in terms of the relative entropy between $\mathbb{Q}$ and $\mathbb{P}$ as a single constraint on the entire path of perturbation. Then, {by considering} $h$ as an adversarial control process in a two-player zero-sum stochastic differential game (in the sense of Elliot and Kalton \cite{EllK72}), {we can reformulate the optimization problem in} \eqref{Eq2.11} as an ``{\it min-max}" problem with uncertainty-aversion between the {\it follower} and that of the adversary (e.g., see \cite{HanS01} or \cite{CheE02} for related discussions, but in a different context).
\end{remark}

\end{document}